\documentclass[french,a4paper,12pt]{article}
\usepackage[T1]{fontenc}
\usepackage[english,frenchb]{babel}
\usepackage{amsmath,amsthm}
\usepackage{amssymb}
\usepackage[all]{xy}
\usepackage{stmaryrd}
\usepackage[francais]{layout}
\theoremstyle{definition}
\newtheorem{theo}{Theorem}

\theoremstyle{definition}

\newtheorem{cor}[theo]{Corollary}
\newtheorem{prop}[theo]{Proposition}

\newcommand{\ind}{\mathop{\mathrm{ind}}}
\newcommand{\Exp}{\mathop{\mathrm{exp}}}
\newcommand{\Br}{\mathop{\mathrm{Br}}}
\newcommand{\Spec}{\operatorname{Spec}}
\newcommand{\Mor}{\operatorname{Mor}}
\newcommand{\PGL}{\operatorname{\mathrm{PGL}}}
\newcommand{\CM}{\operatorname{CM}}

\author{\small{Charles De Clercq}}
\title{\large{\textbf{Classification of upper motives of algebraic groups of\\ inner type $A_n$\\
Classification des motifs supérieurs des groupes algébriques intérieurs de type $A_n$.}}}
\date{}

\begin{document}
\maketitle
\abstract{Soient $A$, $A'$ deux algèbres centrales simples sur un corps $F$ et $\mathbb{F}$ un corps fini de caractéristique $p$. Nous prouvons que les facteurs directs indécomposables supérieurs des motifs de deux variétés anisotropes de drapeaux d'idéaux à droite $X(d_1,...,d_k;A)$ et $X(d'_1,...,d'_s;A')$ à coefficients dans $\mathbb{F}$ sont isomorphes si et seulement si les valuations $p$-adiques de $\mathrm{pgcd}(d_1,...,d_k)$ et $\mathrm{pgcd}(d'_1,..,d'_s)$ sont égales et les classes des composantes $p$-primaires $A_p$ et $A'_p$ de $A$ et $A'$ engendrent le même sous-groupe dans le groupe de Brauer de $F$. Ce résultat mène à une surprenante dichotomie entre les motifs supérieurs des groupes algébriques absolument simples, adjoints et intérieurs de type $A_n$.\\
\begin{center}
\textbf{Abstract}
\end{center}
Let $A$, $A'$ be two central simple algebras over a field $F$ and $\mathbb{F}$ be a finite field of characteristic $p$. We prove that the upper indecomposable direct summands of the motives of two anisotropic varieties of flags of right ideals  $X(d_1,...,d_k;A)$ and $X(d'_1,...,d'_s;A')$ with coefficients in $\mathbb{F}$ are isomorphic if and only if the $p$-adic valuations of $\gcd(d_1,...,d_k)$ and $\gcd(d'_1,..,d'_s)$ are equal and the classes of the $p$-primary components $A_p$ and $A'_p$ of $A$ and $A'$ generate the same group in the Brauer group of $F$. This result leads to a surprising dichotomy between upper motives of absolutely simple adjoint algebraic groups of inner type $A_n$.
}
\selectlanguage{english}
\section{Introduction}
Throughout this note $p$ will be a prime and $\mathbb{F}$ will be a finite field of characteristic $p$. Let $F$ be a field, $F$-$\mathfrak{alg}$ be the category of commutative $F$-algebras and $\CM(F;\mathbb{F})$ be the category of Grothendieck Chow motives with coefficients in $\mathbb{F}$. Motivic properties of projective homogeneous $F$-varieties and their relations with classical discrete invariants have been intensively studied recently (see for example \cite{wit}, \cite{sempetro}, \cite{sempetrzain}, \cite{vish}, \cite{vish2}, \cite{vish3}). In this article we deal with the particular case of projective homogeneous $F$-varieties under the action of an absolutely simple affine adjoint algebraic group of inner type $A_n$. More precisely we prove the following result:

\begin{theo}\label{theo}Let $A$ and $A'$ be two central simple $F$-algebras. The upper indecomposable direct summands of the motives of two anisotropic varieties of flags of right ideals $X(d_1,...,d_k;A)$ and $X(d'_1,...,d'_s;A')$ in $\CM(F;\mathbb{F})$ are isomorphic if and only if $v_p(\gcd(d_1,...,d_k))=v_p(\gcd(d'_1,..,d'_s))$ and the $p$-primary components $A_p$ and $A'_p$ of $A$ and $A'$ generate the same subgroup of $\Br(F)$.
\end{theo}
In $\S 1$ we recall classical discrete invariants and varieties associated to central simple $F$-algebras, while $\S 2$ is devoted
to the theory of upper motives. Finally we prove theorem \ref{theo} in $\S 3$, using an index reduction formula due to Merkurjev, Panin and Wadsworth and the theory of upper motives. Theorem \ref{theo} gives a purely algebraic criterion to compare upper direct summands of varieties of flags of ideals, and leads to a quite unexpected dichotomy between upper motives of absolutely simple adjoint algebraic groups of inner type $A_n$.

\section{Generalities on central simple algebras}

Our reference for classical notions on central simple $F$-algebras is \cite{KMRT}. A finite-dimensional $F$-algebra $A$ is a central simple $F$-algebra if its center $Z(A)$ is equal to $F$ and if $A$ has no non-trivial two-sided ideals. The square root of the $F$-dimension of $A$ is the \emph{degree} of $A$, denoted by $\deg(A)$. Two central simple $F$-algebras $A$ and $B$ are \emph{Brauer-equivalent} if $M_n(A)$ and $M_m(B)$ are isomorphic for some integers $n$ and $m$, and the \emph{Schur index} $\ind(A)$ of a central simple $F$-algebra $A$ is the degree of the (uniquely determined up to isomorphism) central division $F$-algebra Brauer-equivalent to $A$. The tensor product endows the set $\Br(F)$ of equivalence classes of central simple $F$-algebras under the Brauer equivalence with a structure of a torsion abelian group. The exponent of $A$, denoted by $\Exp(A)$, is the order of the class of $A$ in $\Br(F)$ and divides $\ind(A)$.

Let $A$ be a central simple $F$-algebra and $0\leq d_1<...<d_k\leq\deg(A)$ be a sequence of integers. A convenient way to define the variety of flags of right ideals of reduced dimension $d_1$,..., $d_k$ in $A$ is to use the language of functor of points. For any $R$ in $F$-$\mathfrak{alg}$, the set of $R$-points $\Mor_F\left(\Spec(R),X(d_1,...,d_k;A)\right)$ consists of the sequences $(I_1,...,I_k)$ of right ideals of the Azumaya $R$-algebra $A\otimes_F R$ such that $I_1\subset...\subset I_k$, the injection of $A_R$ modules $I_s\rightarrow A_R$ splits and the rank of the $R$-module $I_s$ is equal to $d_s\cdot \deg(A)$ for any $1\leq s \leq k$. For any morphism $R\rightarrow S$ of $F$-algebras the corresponding map from $R$-points to $S$-points is given by $(I_1,...,I_k)\mapsto (I_1\otimes_R S,...,I_k\otimes_R S)$. Two important particular cases of varieties of flags of right ideals are the classical Severi-Brauer variety $X(1;A)$, and the generalized Severi Brauer varieties $X(i;A)$. If $G$ is an algebraic group and $X$ a projective $G$-homogeneous $F$-variety, we say that $X$ is \emph{isotropic} if $X$ has a zero-cycle of degree coprime to $p$, and $X$ is \emph{anisotropic} if $X$ is not isotropic. If $X=X(d_1,...,d_k;A)$ is a variety of flags of right ideals, $X$ is isotropic if and only if $v_p(\gcd(d_1,...,d_k))\geq v_p(\ind(A))$. Note that if the Schur index of $A$ is a power of $p$, $X$ is isotropic if and only if $X$ has a rational point.

\section{The theory upper motives}

Our basic references for the definitions and the main properties of Chow groups with coefficients and the category $\CM(F;\Lambda)$ of Grothendieck Chow motives with coefficients in a commutative ring $\Lambda$ are \cite{andre} and \cite{EKM}. In the sequel $G$ will be a semisimple affine adjoint algebraic group of inner type, $X$ a projective $G$-homogeneous $F$-variety and $\Lambda$ will be assumed to be a finite and connected ring. By \cite{chermer} (see also \cite{upper}) the motive of $X$ decomposes in a unique way (up to isomorphism) as a direct sum of indecomposable motives under these assumptions. Among all the indecomposable direct summands in the complete motivic decomposition of $X$, the (uniquely determined up to isomorphism) indecomposable direct summand $M$ such that the $0$-codimensional Chow group of $M$ is non-zero is the \emph{upper motive} of $X$.

Upper motives are essential : any indecomposable direct summand in the complete motivic decomposition of $X$ is the upper motive of another projective $G$-homogeneous $F$-variety by \cite[Theorem 3.5]{upper}. This structural result implies that the study of the motivic decomposition of a projective $G$-homogeneous $F$-variety $X$ is reduced to the case $\Lambda=\mathbb{F}_p$. Indeed by \cite[Corollary 2.6]{vish4} the complete motivic decomposition of $X$ with coefficients in $\Lambda$ remains the same when passing to the residue field of $\Lambda$, and is also the same as if the ring of coefficients is $\mathbb{F}_p$ by \cite[Theorem 2.1]{decvarprojcoeff}, where $p$ is the characteristic of the residue field of $\Lambda$. These results motivate the study of the set $\mathfrak{X}_G$ of \emph{upper $p$-motives} of the algebraic group $G$, which consists of the isomorphism classes of upper motives of projective $G$-homogeneous $F$-varieties in $\CM(F;\mathbb{F}_p)$. Furthermore the dimension of the upper motive of $X$ in $\CM(F;\mathbb{F}_p)$ is equal to the canonical $p$-dimension of $X$ by \cite[Theorem 5.1]{ICM}, hence upper motives encode important information on the underlying varieties. Upper motives also have good properties : the upper motives of two projective $G$-homogeneous $F$-varieties $X$ and $X'$ in $\CM(F;\mathbb{F})$ are isomorphic if and only if both $X_{F(X')}$ and $X'_{F(X)}$ are isotropic by \cite[Corollary 2.15]{upper}. The variety $X$ is isotropic if and only if the upper motive of $X$ is isomorphic to the \emph{Tate motive} (that is to say the motive of $\Spec(F)$) and this is why we focus in this note on the case of anisotropic varieties of flags of right ideals. 

If $G$ is absolutely simple adjoint of inner type $A_n$, $G$ is isomorphic to $\PGL_1(A)$, where $A$ is a central simple $F$-algebra of degree $n+1$. Any projective $G$-homogeneous $F$-variety is then isomorphic to a variety $X(d_1,...,d_k;A)$ of flags of right ideals in $A$ (see \cite{index}) thus theorem \ref{theo} classifies upper motives of absolutely simple affine adjoint algebraic groups of inner type $A_n$. In the particular case of classical Severi-Brauer varieties theorem \ref{theo} corresponds to \cite[Theorem 9.3]{amitsur}, since for any field extension $E/F$ a central simple $F$-algebra split over $E$ if and only if the Severi-Brauer variety $SB(1,A_E)$ has a rational point.

\section{Main results}

Let $D$ be a central division $F$-algebra of degree $p^n$. For any $0\leq k\leq n$, $M_{k,D}$ will denote the upper indecomposable direct summand of $X(p^k;D)$ in $\CM(F;\mathbb{F})$. If $D'$ is another central division $F$-algebra of degree $p^n$ and $j$ satisfies $1\leq j \leq p^n$, we denote the integer $\frac{p^k}{\gcd(j,p^k)}\cdot \ind(D\otimes D'^{-j})$ by $\mu_{k,j}^{D,D'}$. In the sequel the following index reduction formula (see \cite[p. 565]{index})  will be of constant use :
$$\ind(D_{F(X(p^k;D'))})=\gcd_{1\leq j\leq p^n}\mu_{k,j}^{D,D'}=\min_{1\leq j\leq p^n}\mu_{k,j}^{D,D'}$$ 

\begin{prop}\label{prop1}
Let $D$ and $D'$ be two central division $F$-algebras of degree $p^n$. Assume that $\Exp(D)\geq \Exp(D')$ and that $X(p^k;D)_{F(X(p^k;D'))}$ is isotropic for some integer $0\leq k<n$. If $\ind(D_{F(X(k;D'))})=\mu_{k,j_0}^{D,D'}$, $j_0$ is coprime to $p$.
\end{prop}
\begin{proof}
Suppose that $p$ divides $j_0$ and $\ind(D_{F(X(k;D'))})=\mu_{k,j_0}^{D,D'}$. By assumption $X(k;D)_{F(X(k;D'))}$ has a rational point, hence the integer $\mu_{k,j_0}^{D,D'}$ divides $p^k$ by \cite[Proposition 1.17]{KMRT} and $\ind(D\otimes D'^{-j_0})\mid \gcd(j_0,p^k)$. Since $p$ divides $j_0$, $\Exp(D'^{-j_0})<\Exp(D')$, therefore $\Exp(D'^{-j_0})<\Exp(D)$ and $\Exp(D)=\Exp(D\otimes D'^{-j_0})$. It follows that $\Exp(D)$ divides $j_0$, thus $\Exp(D')$ also divides $j_0$. The central simple $F$-algebra $D'^{j_0}$ is therefore split and $D\otimes D'^{-j_0}$ is Brauer-equivalent to $D$ so that $\ind(D)$ divides $p^k$, a contradiction.
\end{proof}

\begin{theo}\label{theo1}
Let $\mathbb{F}$ be a finite field of characteristic $p$ and $D$, $D'$ be two central division $F$-algebras of degree $p^n$. The following assertions are equivalent :\\
1) for some integer $0\leq l<n$, $M_{l,D}$ and $M_{l,D'}$ are isomorphic in $\CM(F;\mathbb{F})$;\\
2) the classes of $D$ and $D'$ generate the same subgroup of $\Br(F)$;\\
3) for any $0\leq l<n$, $M_{l,D}$ is isomorphic to $M_{l,D'}$ in $\CM(F;\mathbb{F})$.
\end{theo}

\begin{proof}We first show that $1)\Rightarrow 2)$. We may replace $D$ by $D'$ and thus assume that $\Exp(D)$ is greater than $\Exp(D')$. Since $M_{l,D}$ is isomorphic to $M_{l,D'}$, there is an integer $j_0$ coprime to $p$ such that the Schur index of $D\otimes D'^{-j_0}$ is equal to $1$ by \cite[Proposition 1.17]{KMRT} and proposition \ref{prop1}, hence $D\otimes D'^{-j_0}$ is split and the class of $D$ is equal to the class of $D'^{j_0}$ in $\Br(F)$. Furthermore since $j_0$ is coprime to $p$ the class of $D$ in $\Br(F)$ is also a generator of the subgroup of $\Br(F)$ generated by $[D']$. Now statement $2)\Rightarrow 3)$ : if $[D]$ and $[D']$ generate the same group in $\Br(F)$, $\ind(D_E)=\ind(D'_E)$ for any field extension $E/F$. Given an integer $0\leq l<n$, since $X(p^l;D)$ has a rational point over $F(X(p^l;D))$, $\ind(D'_{F(X(p^l;D))})=\ind(D_{F(X(p^l;D))})$ divides $p^l$. The variety $X(p^l;D')$ then also has a rational point over $F(X(p^l;D))$ by \cite[Proposition 1.17]{KMRT}. The same procedure replacing $D$ by $D'$ shows that $X(p^l;D)$ has a rational point over $F(X(p^l;D'))$, hence $M_{l,D}$ is isomorphic to $M_{l,D'}$. Finally $3)\Rightarrow 1)$ is obvious.
\end{proof}

\begin{cor}\label{cor1}
Let $D$ and $D'$ be two central division $F$-algebras of degree $p^n$ and $p^{n'}$. The upper summands $M_{k,D}$ and $M_{l,D'}$ are isomorphic for some integers $0\leq k<n$ and $0\leq l<n'$ if and only if $k=l$ and the classes of $D$ and $D'$ generate the same subgroup of $\Br(F)$.
\end{cor}
\begin{proof}
Since by \cite[Theorem 4.1]{upper} the generalized Severi-Brauer varieties $X(p^k;D)$ and $X(p^l;D')$ are $p$-incompressible, if $M_{k,D}$ and $M_{l,D'}$ are isomorphic, the dimension of $X(p^k;D)$ (which is $p^k(p^n-p^k)$) is equal to the dimension of $X(p^l;D')$. The equality $p^k(p^n-p^k)=p^l(p^{n'}-p^l)$ implies that $k=l$, $n=n'$ and it remains to apply theorem \ref{theo1}. The converse is clear by theorem \ref{theo1}. 
\end{proof}
\begin{proof}[Proof of theorem \ref{theo}] Set $X=X(d_1,...,d_k;A)$, $Y=X(d_1,...,d_k;A')$, and also $u=v_p(\gcd(d_1,...,d_k))$ and $v=v_p(\gcd(d'_1,..,d'_s))$. If $D$ and $D'$ are two central division $F$-algebras Brauer-equivalent to $A_p$ and $A'_p$, the upper indecomposable direct summand of $X$ (resp. of $Y$) is isomorphic to $M_{u,D}$ (resp. to $M_{v,D'}$) by \cite[Theorem 3.8]{upper}. By corollary \ref{cor1} these summands are isomorphic if and only if $u=v$ (since $X$ and $Y$ are anisotropic) and the classes of $A_p$ and $A'_p$ generate the same subgroup of $\Br(F)$.
\end{proof}

\begin{theo}\label{theo3}Let $G$ and $G'$ be two absolutely simple affine adjoint algebraic groups of inner type $A_n$ and $A_{n'}$. Then either $\mathfrak{X}_G\cap \mathfrak{X}_{G'}$ is reduced to the class of the Tate motive or $\mathfrak{X}_G=\mathfrak{X}_{G'}$.
\end{theo}
\begin{proof}If $\mathfrak{X}_{\PGL_1(A)}\cap \mathfrak{X}_{\PGL_1(A')}$ is not reduced to the class of the Tate motive, there are two anisotropic varieties of flags of right ideals $X=X(d_1,...,d_k;A)$ and $Y=X(d'_1,...,d'_s;A')$ whose upper motives are isomorphic. By theorem \ref{theo} this implies that the upper $p$-motive of any anisotropic $\PGL_1(A)$-homogeneous $F$-variety $X(d_1,...,d_s;A)$ is isomorphic to, say, the upper $p$-motive of $X(d_1,...,d_s;A')$.
\end{proof}

\textbf{Aknowledgements :} I would like to express my gratitude to N. Karpenko, for introducing me to this subject, raising this question and for stimulating discussions on this subject.
\begin{footnotesize}

\end{footnotesize}
\end{document}